\centerline{\bf Computational solutions of unified fractional reaction-diffusion equations}
\vskip.1cm\centerline{\bf with composite fractional time derivative}

\vskip.3cm \centerline{ \bf R.K. Saxena} 
\vskip0cm\centerline{Department of Mathematics and Statistics} 
\vskip0cm\centerline{Jai Narain Vyas University, Jodhpur-342005, India, ram.saxena@yahoo.com}

\vskip.3cm\centerline{\bf A.M. Mathai} 
\vskip0cm\centerline{Centre for Mathematical Sciences, Arunapuram P.O., Pala,}
\vskip0cm\centerline{Kerala-686 574, India, directorcms458@gmail.com}
\vskip0cm\centerline{and} 
\vskip0cm\centerline{Department of Mathematics and Statistics, McGill University,}
\vskip0cm\centerline{Montreal, Canada H3A 2K6, mathai@math.mcgill.ca} 

\vskip.3cm\centerline{\bf H.J. Haubold} 
\vskip0cm\centerline{Office for Outer Space Affairs, United Nations, Vienna International Centre,}
\vskip0cm\centerline{P.O. Box 500, A-1400 Vienna, Austria, hans.haubold@unvienna.org}
\vskip0cm\centerline{and} 
\vskip0cm\centerline{Centre for Mathematical Sciences, Arunapuram P.O., Pala,}
\vskip0cm\centerline{Kerala-686 574, India}

\vskip.5cm\noindent{\bf Abstract} \vskip.3cm This paper deals with the
investigation of the computational solutions of an unified
fractional reaction-diffusion equation, which is obtained from the
standard diffusion equation by replacing the time derivative of
first order by the generalized fractional time-derivative defined by
Hilfer (2000), the space derivative of second order by the
Riesz-Feller fractional derivative and adding the function
$\phi(x,t)$ which is a nonlinear function governing reaction. 
The solution is derived by the application of
the Laplace and Fourier transforms in a compact and closed form in
terms of the H-function. The main result obtained in this paper
provides an elegant extension of the fundamental solution for the
space-time fractional diffusion equation obtained earlier by
Mainardi et al. (2001, 2005) and a result very recently given by
Tomovski et al. (2011). Computational representation of the
fundamental solution is also obtained explicitly. Fractional order
moments of the distribution are deduced. At the end, mild extensions
of the derived results associated with a finite number of
Riesz-Feller space fractional derivatives are also discussed.



\vskip.5cm\noindent{\bf 1.\hskip.3cm Introduction} \vskip.3cm
Standard (not fractional) reaction-diffusion equations are an
important class of partial differential equations to investigate
nonlinear behavior of complex systems. Standard nonlinear reaction-diffusion equations
can be simulated by numerical techniques, such as finite difference
methods. The reaction-diffusion equation takes into account particle
diffusion (diffusion constant and spatial Laplacian operator) and
particle reaction (reaction constants and nonlinear reactive terms).
Well-known special cases of such standard reaction-diffusion
equations are the (i) Schloegl model, (ii) Fisher-Kolmogorov
equation, (iii) real and complex Ginzburg-Landau equations, (iv)
FitzHugh-Nagumo model, and (v) Gray-Scott model. These equations are
known under their respective names both in the mathematical and
physical literature. The nontrivial behavior of these equations
arises from the competition between the reaction/relaxation and
diffusion/transport. In recent years, interest is developed by several authors
in the applications of reaction-diffusion models in pattern
formation in physical and biological sciences and for describing 
non-Gaussian, non-Markovian, and non-Fickian phenomena in complex
systems. In this connection,
one can refer to Murray (2003), Kuramoto (2003), Wilhelmsson and
Lazzaro (2001), and Hundsdorfer and Verwer (2003). These systems
show that diffusion can produce the spontaneous formation of
spatio-temporal patterns. For details, see the work of Nicolis and
Prigogine (1977) and Haken (2004). A general model for
reaction-diffusion systems is investigated by Henry and Wearne
(2000, 2002) and Henry et al. (2005), and Haubold et al. (2007, 2011, 2012).

\vskip.2cm In this paper, we investigate the solution of an unified
model of reaction-diffusion system (3.1) in which the two-parameter
fractional derivative ${_0D}_{0_{+}}^{\mu,\nu}$ acts as a
time-derivative and the Riesz-Feller derivative
${_xD}_{\theta}^{\alpha}$ as the space-derivative. This  new model
provides an extension of the models discussed earlier by Jesperson
et al. (1999), Del-Castillo-Negrete et al. (2003), Mainardi et al.
(2001, 2005), Kilbas et al. (2004), Haubold et al. (2007),
Saxena (2012), Saxena et al. (2010), Saxena et al. (2006
a,b,c,d), and Tomovski et al. (2001). The results are obtained in a
compact form, which are suitable for numerical computation.
Computational representation of the fundamental solution is
explicitly derived and fractional order moments are investigated.
For recent and related works on fractional kinetic equations and
reaction-diffusion problems, one can refer to papers by Haubold and
Mathai (1995, 2000) and Saxena et al. (2002, 2004a,b,c, 2006a,d,
2008).

\vskip.3cm\noindent{\bf 2.\hskip.3cm Unified fractional
reaction-diffusion equations}
\vskip.3cm In this section, we will
investigate the solution of the unified reaction-diffusion model
(2.1). The main result is given in the form of the following

\vskip.3cm\noindent{\bf Theorem 2.1}.\hskip.3cm{\it Consider an
unified fractional reaction-diffusion model
$${_0D}_t^{\mu,\nu}N(x,t)=\eta~{_xD}_{\theta}^{\alpha}N(x,t)+\phi(x,t)\eqno(2.1)
$$where $\eta,t>0,x\in R; \alpha,\theta,\beta$ are real parameters
with the constraints  $0<\alpha\le \min (\alpha,2-\alpha),$
${_0D}_t^{\mu,\nu}$ is the generalized Riemann-Liouville fractional
derivative operator defined by (A7) with the initial conditions
$$[I_{0_{+}}^{(1-\nu)(1-\mu),0}]N(x,0_{+})=N_0(x);~0<\mu<1,~0\le
\nu\le 1\eqno(2.2)
$$involving the Riemann-Liouville fractional integral of order
$(1-\nu)(1-\mu)$ evaluated for $t\to 0_{+}$ and  $\lim_{|x|\to
\infty}N(x,t)=0.$ Here ${_xD}_{\theta}^{\alpha}$ is the Riesz-Feller
space-fractional derivative of order $\alpha$ and asymmetry $\theta$
defined by (A11), $\eta$ is a diffusion constant and $\phi(x,t)$ is
a nonlinear function belonging to the area of reaction-diffusion.
Then for the solution of (2.1) subject to the above constraints,
there holds the formula
$$\eqalignno{N(x,t)&=\int_0^xG_1(x-\tau,t)N_0(\tau){\rm d}\tau\cr
&+\int_0^t(t-\xi)^{\beta-1}\int_0^{\infty}G_2(x-\tau,t-\xi)\phi(\tau,\xi){\rm
d}\tau~{\rm d}\xi&(2.3)\cr \noalign{\hbox{where}}
G_1(x,t)&={{t^{\mu+\nu(1-\mu)-1}}\over{\alpha|x|}}H_{3,3}^{2,1}\left[{{|x|}\over{\eta^{{1}\over{\alpha}}t^{{\mu}
\over{\alpha}}}}\bigg\vert_{(1,{{1}\over{\alpha}}),(1,1),(1,\rho)}^{(1,{{1}\over{\alpha}}),(\mu+\nu(1-\mu),{{\mu}\over{\alpha}}),(1,\rho)}\right]&(2.4)\cr
\noalign{\hbox{$\alpha>0$ and }}
G_2(x,t)&={{1}\over{\alpha|x|}}H_{3,3}^{2,1}\left[{{|x|}\over{\eta^{{1}\over{\alpha}}t^{{\mu}\over{\alpha}}}}
\bigg\vert_{(1,{{1}\over{\alpha}}),(1,1),(1,\rho)}^{(1,{{1}\over{\alpha}}),(\mu,{{\mu}\over{\alpha}}),(1,\rho)}\right]&(2.5)\cr}
$$for $\alpha>0$ where $H_{3,3}^{2,1}$ is the familiar H-function
[Mathai and Saxena (1978), Kilbas and Saigo (2004) and Mathai et al.
(2010)], $\Re(\mu)>0,\Re(\mu+\nu(1-\mu))>0.$}

\vskip.3cm\noindent{\bf Proof:}\hskip.3cm If we apply the Laplace
transform (Erd\'elyi, et al. 1954) with respect to the time variable
$t$ and Fourier transform with respect to space variable $x$ and use
the initial conditions  and the formulas (A8) and (A11), then the
given equation transforms into the form
$$s^{\mu}\tilde{N}^{*}(k,s)-s^{\nu(\mu-1)}N^{*}_0(k)=-\eta~\psi_{\alpha}^{\theta}(k)\tilde{N}^{*}(k,s)+\phi^{*}(k,s),\eqno(2.6)
$$where according to the conventions followed, the symbol $\sim$
will stand for the Laplace transform with respect to time variable
$t$ and * represents the Fourier transform with respect to space
variable $x$. Solving for $\tilde{N}^{*}(k,s)$, it yields
$$\tilde{N}^{*}(k,s)={{N_0^{*}(k)s^{\nu(\mu-1)}}\over{s^{\mu}+\eta~\psi_{\alpha}^{\theta}(k)}}+{{\tilde{\phi}^{*}(k)}
\over{s^{\mu}+\eta~\psi_{\alpha}^{\theta}(k)}}.\eqno(2.7)
$$On taking the inverse Laplace transform of (2.7) and applying the
formula
$$\eqalignno{L^{-1}\{{{s^{\beta-1}}\over{a+s^{\alpha}}}\}&=t^{\alpha-\beta}E_{\alpha,\alpha-\beta+1}(-at^{\alpha}),&(2.8)\cr
\noalign{\hbox{where $\Re(s)>0,\Re(\alpha)>0,\Re(\alpha-\beta)>-1$,
it is seen that}}
N^{*}(k,t)&=N_0^{*}(k)t^{\mu+\nu(1-\mu)-1}E_{\mu,\mu+\nu(1-\mu)}(-\eta~t^{\mu}\psi_{\alpha}^{\theta}(k))\cr
&+\int_0^t\phi^{*}(k,t-\xi)\xi^{\mu-1}E_{\mu,\mu}(-\eta~\psi_{\alpha}^{\theta}(k)\xi^{\mu}){\rm
d}\xi.&(2.9)\cr \noalign{\hbox{Taking the inverse Fourier transform
of (2.9), we find}}
N(x,t)&={{t^{\mu+\nu(1-\mu)-1}}\over{2\pi}}\int_{-\infty}^{\infty}N_0^{*}(k)E_{\mu,\mu+\nu(1-\mu)}(-\eta~t^{\mu}\psi_{\alpha}^{\theta}(k))\exp(-ikx){\rm
d}k\cr
&+{{1}\over{2\pi}}\int_0^t\xi^{\mu-1}\int_{-\infty}^{\infty}\phi^{*}(k,t-\xi)E_{\mu,\mu}(-\eta
t^{\mu}\psi_{\alpha}^{\theta}(k))\exp(-ikx){\rm d}k~{\rm
d}\xi.&(2.10)\cr}
$$If we now apply the convolution theorem of the Fourier transform
to (2.10) and make use of the following inverse Fourier transform
formula (Haubold et al. 2007):
$$F^{-1}[E_{\beta,\gamma}(-\eta
t^{\beta}\psi_{\alpha}^{\theta}(k);x]={{1}\over{\alpha
|x|}}H_{3,3}^{2,1}\left[{{|x|}\over{\eta^{{1}\over{\alpha}}t^{{\beta}\over{\alpha}}}}
\bigg\vert_{(1,{{1}\over{\alpha}}),(1,1),(1,\rho)}^{(1,{{1}\over{\alpha}}),(\gamma,{{\beta}\over{\alpha}}),(1,\rho)}\right],\eqno(2.11)
$$where $\Re(\alpha)>0,\Re(\beta)>0,\Re(\gamma)>0$ and $\rho
=({{\alpha-\theta}\over{2\alpha}})$, it gives the desired solution
in the form
$$\eqalignno{N(x,t)&=\int_0^xG_1(x-\tau,t)N_0(\tau){\rm d}\tau\cr
&+\int_0^t(t-\xi)^{\beta-1}\int_0^xG_2(x-\tau,t-\xi)\phi(\tau,\xi){\rm
d}\tau~{\rm d}\xi&(2.12)\cr \noalign{\hbox{where}}
G_1(x,t)&={{t^{\mu+\nu(1-\mu)-1}}\over{2\pi}}\int_{-\infty}^{\infty}\exp(-ikx)E_{\mu,\mu+\nu(1-\mu)}(-\eta
t^{\mu}\psi_{\alpha}^{\theta}(k)){\rm d}k\cr
&={{t^{\mu+\nu(1-\mu)-1}}\over{\alpha|x|}}H_{3,3}^{2,1}\left[{{|x|}\over{\eta^{{1}\over{\alpha}}t^{{\mu}\over{\alpha}}}}
\bigg\vert_{(1,{{1}\over{\alpha}}),(1,1),(1,\rho)}^{(1,{{1}\over{\alpha}}),(\mu+\nu(1-\mu),
{{\mu}\over{\alpha}}),(1,\rho)}\right]&(2.13)\cr \noalign{\hbox{for
$\alpha >0$ and }}
G_2(x,t)&={{1}\over{2\pi}}\int_{-\infty}^{\infty}\exp(-ikx)E_{\mu,\mu}(-\eta
t^{\mu}\psi_{\alpha}^{\theta}(k)){\rm d}k\cr
&={{1}\over{\alpha|x|}}H_{3.3}^{2,1}\left[{{|x|}\over{\eta^{{1}\over{\alpha}}t^{{\mu}\over{\alpha}}}}
\bigg\vert_{(1,{{1}\over{\alpha}}),(1,1),(1,\rho)}^{(1,{{1}\over{\alpha}}),(\mu,{{\mu}\over{\alpha}}),(1,\rho)}\right],&(2.14)\cr}
$$where $\Re(\alpha)>0,\Re(\mu)>0,\Re(\mu+\nu(1-\mu))>0$. This
completes the proof of the theorem.

\vskip.3cm\noindent{\bf 3.\hskip.3cm Special cases}

\vskip.3cm When $\theta=0$, the Riesz-Feller space derivative
reduces to the Riesz space fractional derivative, and consequently
we arrive at the following result recently given by Tomovski et al.
(2011).

\vskip.3cm\noindent{\bf Corollary 3.1}.\hskip.3cm The solution of
extended fractional reaction-diffusion equation
$$\eqalignno{{_0D}_{0_{+}}^{\mu,\nu}N(x,t)-\eta~{_xD}_0^{\alpha}N(x,t)&=\phi(x,t),~x\in
R,~t>0,~\eta>0&(3.1)\cr \noalign{\hbox{with initial conditions}}
0\le \nu\le 1,~0<\mu<1\hbox{ and }&0<\alpha\le
2,~\lim_{x\to\pm\infty}N(x,t)=0&(3.2)\cr
[I_{0_{+}}^{(1-\nu)(\mu-1),0}]N(x,0_{+})&=N_0(x)&(3.3)\cr}
$$where $\eta$ is a diffusion constant, ${_0D}_{0_{+}}^{\mu,\nu}$ is
the generalized Riemann-Liouville fractional derivative operator
defined by (A7), ${_xD_0}^{\alpha}$ is the Riesz space fractional
derivative operator defined by (A15) and $\phi(x,t)$ is a nonlinear
function belonging to the area of reaction-diffusion, is given by
$$\eqalignno{N(x,t)&=\int_0^xG_3(x-\tau,t)N_0(\tau){\rm d}\tau\cr
&+\int_0^t(t-\xi)^{\beta-1}\int_0^xG_4(x-\tau,t-\xi)\phi(\tau,\xi){\rm
d}\tau~{\rm d}\xi,&(3.4)\cr \noalign{\hbox{where}}
G_3(x,t)&={{t^{\mu+\nu(1-\mu)-1}}\over{\alpha|x|}}H_{3,3}^{2,1}\left[{{|x|}\over{\eta^{{1}\over{\alpha}}t^{{\mu}\over{\alpha}}}}
\bigg\vert_{(1,{{1}\over{\alpha}}),(1,1),(1,{1\over2})}^{(1,{{1}\over{\alpha}}),(\mu+\nu(1-\mu),{{\mu}\over{\alpha}}),(1,{1\over2})}\right]&(3.5)\cr
\noalign{\hbox{for $\alpha>0$ and }}
G_4(x,t)&={{1}\over{\alpha|x|}}H_{3.3}^{2.1}\left[{{|x|}\over{\eta^{{1}\over{\alpha}}t^{{\mu}\over{\alpha}}}}
\bigg\vert_{(1,{{1}\over{\alpha}}),(1,1),(1,{1\over2})}^{(1,{{1}\over{\alpha}}),(\mu,{{\mu}\over{\alpha}}),(1,{1\over2})}\right],~\alpha>0.&(3.6)\cr}
$$
 \vskip.3cm\noindent{\bf Note 3.1.}\hskip.3cm Expressions for $G_3$
and $G_4$ can be obtained from $G_1$ and $G_2$ respectively by
taking $\rho={1\over2}$. \vskip.2cm Similarly, if we set
$N_0(x)=\delta(x), \phi=0$, where $\delta (x)$ is the Dirac delta
function, then the theorem reduces to the following interesting
result: \vskip.3cm\noindent{\bf Corollary 3.2.}\hskip.3cm {\it
Consider the following reaction-diffusion model

$$\eqalignno{{_0D}_{0_{+}}^{\mu,\nu}N(x,t)&=\eta~{_xD}_{\theta}^{\alpha}N(x,t),~\eta>0,~x\in
R, ~0<\alpha\le 2,&(3.7)\cr \noalign{\hbox{with the initial
conditions}}
(I_{0_{+}}^{(1-\nu)(\mu-1),0})N(x,0_{+})&=\delta(x),~0<\mu<1,~0\le
\nu\le 1,~\lim_{x\to\pm\infty}N(x,t)=0&(3.8)\cr}
$$where $\eta$ is a diffusion constant and $\delta(x)$ is the Dirac
delta function. Then for the fundamental solution of (3.7) with the
initial conditions (3.8), there holds the formula
$$N(x,t)={{t^{\mu+\nu(1-\mu)-1}}\over{\alpha|x|}}H_{3,3}^{2,1}\left[{{|x|}\over{\eta^{{1}\over{\alpha}}t^{{\mu}\over{\alpha}}}}
\bigg\vert_{(1,{{1}\over{\alpha}}),(1,1),(1,\rho)}^{(1,{{1}\over{\alpha}}),(\mu+\nu(1-\mu),{{\mu}\over{\alpha}}),(1,\rho)}\right]\eqno(3.9)
$$for $\alpha>0$ where $\rho={{\alpha-\theta}\over{2\alpha}}$.}

\vskip.3cm Some interesting special cases of  (3.9) are enumerated
below. \vskip.2cm\noindent (i): When $\nu=1$, then the operator
$D_{0_{+}}^{\mu,\nu}$ reduces to Caputo fractional derivative
$^C_0D_t^{\mu}$ defined in (A5) and the result (3.9) yields the
fundamental solution of the space-time fractional diffusion equation
obtained by Mainardi et al. (2001) and
Mainardi et al. (2005). \vskip.2cm\noindent(ii): {\bf
Neutral fractional diffusion}. We note that for $\mu=\alpha$, the
solution (3.9) reduces to
$$N_{\alpha,\nu}^{\theta}(x,t)={{t^{\alpha+\nu(1-\alpha)-1}}\over{\alpha|x|}}H_{3,3}^{2,1}\left[{{|x|}\over{(t~\eta^{{1}\over{\alpha}})}}
\bigg\vert_{(1,{{1}\over{\alpha}}),(1,1),(1,\rho)}^{(1,{{1}\over{\alpha}}),(\alpha+\nu(1-\alpha),1),(1,\rho)}\right].\eqno(3.10)
$$If we further set $\nu=1$, then the time-fractional operator
becomes the Caputo operator and neutral fractional diffusion occurs.
It  will be denoted by the conventional symbol
$N_{\alpha}^{\theta}$, (3.10) now simplifies into
$$N_{\alpha}^{\theta}(x,t)={{1}\over{\alpha|x|}}H_{2,2}^{1,1}\left[{{|x|}\over{t~\eta^{{1}\over{\alpha}}}}
\bigg\vert_{(1,{{1}\over{\alpha}}),(1,\rho)}^{(1,{{1}\over{\alpha}}),(1,\rho)}\right],\eqno(3.11)
$$which can be expressed in terms of a Mellin-Barnes type integral
as
$$N_{\alpha}^{\theta}(x,t)={{1}\over{\pi\alpha x}}{{1}\over{2\pi
i}}\int_{\gamma-i\infty}^{\gamma+i\infty}\Gamma({{s}\over{\alpha}})\Gamma(1-{{s}\over{\alpha}})
\sin[{{\pi}\over2}{{s}\over{\alpha}}(\alpha-\theta)](x/t)^s{\rm
d}s,\eqno(3.12)
$$ $(\eta=1$, taken for simplicity). If the poles of the gamma
functions occurring in the integrand of (3.12) are all simple, then
evaluating the integral as a sum of the residues at the simple poles
of $\Gamma(s/\alpha)$ at the points $s=-\alpha n,~ n\in N_0$ and
$\Gamma(1-s/\alpha)$ at the points $s=\alpha+\alpha n, n\in N_0$ we
obtain the series representations
$$\eqalignno{N_{\alpha}^{\theta}(x,t)&={{1}\over{\pi
x}}\sum_{n=0}^{\infty}\sin[{{n\pi}\over{2}}(\theta-\alpha)](-x/t)^{\alpha
n},~~0<x<1&(3.13)\cr \noalign{\hbox{and}}
N_{\alpha}^{\theta}(x,t)&={{1}\over{\pi
x}}\sum_{n=0}^{\infty}\sin[{{n\pi}\over{2}}(\theta-\alpha)^{-\alpha
n}],~~1<x<\infty.&(3.14)\cr}
$$Following the procedure adopted in Gorenflo and Mainardi (1997)
and making use of the formula
$$\eqalignno{\sum_{n=1}^{\infty}r^n\sin(n\alpha)&=\Im
[\sum_{n=1}^{\infty}r^n\exp(in\alpha)]=\Im[{{r\exp(i\alpha)}\over{1-r\exp(i\alpha)}}]\cr
&={{r\sin\alpha}\over{1-2r\cos\alpha+r^2}},\hbox{ for
}|r|<1&(3.15)\cr}
$$where $\alpha\in R$ and it yields the interesting result given by
Mainardi et al. (2001, 2005)
$$N_{\alpha}^{\theta}(x,t)={{1}\over{t\pi}}{{y^{\alpha-1}\sin[{{\pi}\over{2}}(\alpha-\theta)]}\over{1+2y^{\alpha}\cos[{{\pi}\over2}(\alpha-\theta)]+y^{2\alpha}}},~~y={{x}\over{t}},0<x<\infty,
0<\alpha\le 2.\eqno(3.16)
$$Next, we derive some stable densities in terms of the H-functions
as special cases of the solution (3.12). \vskip.2cm\noindent(iii):
If we set $\mu=\nu=1, ~0<\alpha<2,~\theta\le
\min\{\alpha,2-\alpha\}$ then (3.7) reduces to space-fractional
diffusion equation, which we denote by $L_{\alpha}^{\theta}(x)$, is
the fundamental solution of the following space-time fractional
diffusion model:
$${{\partial N(x,t)}\over{\partial
t}}=\eta~{_xD}_{\theta}^{\alpha}N(x,t),~\eta>0, x\in R,\eqno(3.17)
$$with the initial conditions
$N(x,t=0)=\delta(x),\lim_{x\to\pm\infty}N(x,t)=0$, where $\eta$ is
the diffusion constant and $\delta(x)$ is the Dirac delta function.
Hence for the fundamental solution of (3.17) there holds the formula
$$L_{\alpha}^{\theta}(x)={{1}\over{\alpha(\eta
t)^{{1}\over{\alpha}}}}H_{2,2}^{1,1}\left[{{(\eta
t)^{{1}\over{\alpha}}}\over{|x|}}\bigg\vert_{(1,{{1}\over{\alpha}}),(1,\rho)}^{(1,1),(1,\rho)}\right],~0<\alpha<1,~|\theta|\le
\alpha,\eqno(3.18)
$$where $\rho={{\alpha-\theta}\over{2\alpha}}$. The density
represented by the above expression is known as $\alpha$-stable
L\'evy density. By virtue of the H-function formula
(Mathai et al., 2010) another form of this density is given
by
$$L_{\alpha}^{\theta}(x)={{1}\over{\alpha(\eta
t)^{{1}\over{\alpha}}}}H_{2,2}^{1,1}\left[{{|x|}\over{(\eta
t)^{{1}\over{\alpha}}}}\bigg\vert_{(0,1),(0,\rho)}^{(0,{{1}\over{\alpha}}),(0,{{1}\over{\alpha}})}\right],
~1<\alpha\le 2,~|\theta|\le 2-\alpha.\eqno(3.19)
$$
\vskip.3cm\noindent{\bf Remark 3.1}.\hskip.3cm A general
representation of all stable distributions in terms of special
functions has been given by Schneider (1986). It is shown that the
stable probability functions can be expressed by means of the
H-functions, also see Uchaikin and Zolotarev (1999). We further note
that Feller (1952) had derived the representations of the stable
probability functions in terms of convergent and asymptotic power
series in 1952. His results are revisited by Schneider, who had
shown in case of $L_{\alpha}^{\theta}(x)$ to restrict our attention
to $x>0$, since the evaluation for $x<0$ can be done by using the
symmetry property $L_{\alpha}^{\theta}(-x)=L_{\alpha}^{-\theta}(x)$.

\vskip.2cm\noindent(iv):~~Next, if we take
$\alpha=2,~0<\mu<1,~\nu=1,~\theta=0$ then we obtain the
time-fractional diffusion, which is governed by the following time
fractional diffusion model:
$${{\partial^{\mu}}\over{\partial
t^{\mu}}}N(x,t)=\eta~{{\partial^2}\over{\partial
x^2}}N(x,t),~\eta>0,~x\in R,\eqno(3.20)
$$with the initial conditions $N(x,t=0)=\delta(x),
\lim_{x\to\pm\infty}N(x,t)=0$ where $\eta$ is a diffusion constant
and $\delta(x)$ is the Dirac delta function, whose fundamental
solution is given by the equation
$$N(x,t)={{1}\over{2|x|}}H_{1,1}^{1,0}\left[{{|x|}\over{(\eta
t)^{1\over2}}}\bigg\vert_{(1,1)}^{(1,{{\mu}\over2})}\right].\eqno(3.21)
$$(v):~ Further, if we set $\alpha=2,\mu=\nu=1$ and $\theta\to 0$
then for the fundamental solution of the standard diffusion equation
$$\eqalignno{{{\partial}\over{\partial
t}}N(x,t)&=\eta~{{\partial^2}\over{\partial x^2}}N(x,t),&(3.22)\cr
\noalign{\hbox{with initial conditions}}
N(x,t=0)&=\delta(x),~\lim_{x\to\pm\infty}N(x,t)=0,&(3.23)\cr
\noalign{\hbox{there holds the formula}}
N(x,t)&={{1}\over{2|x|}}H_{1,1}^{1,0}\left[{{|x|}\over{(\eta
t)^{1\over2}}}\bigg\vert_{(1,1)}^{(1,{1\over2})}\right]=(4\pi \eta
t)^{-{1\over2}}\exp[-{{|x|^2}\over{4\eta t}}],&(3.24)\cr}
$$which is the classical Gaussian density. For further details and
importance of these special cases based on the Green function, one
can refer to the papers by Mainardi et al. (2001, 2005). For $\nu=0$
the fractional derivative ${_0D}_{0_{+}}^{\mu,\nu}$ reduces to
Riemann-Liouville fractional derivative operator
${^{RL}_0D}_t^{\mu}$, defined by (A3), and the theorem yields

\vskip.3cm\noindent{\bf Corollary 3.3.}\hskip.3cm{\it Consider an
extended fractional reaction-diffusion model
$$\eqalignno{{^{RL}_0D}_t^{\mu}N(x,t)&=\eta~{_xD}_{\theta}^{\alpha}N(x,t)+\phi(x,t)&(3.25)\cr
\noalign{\hbox{where $\eta,t>0,x\in R; \alpha,\theta,\mu$ are real
parameters with the constraints}}
0<\alpha\le 2,&|\theta|\le \min
(\alpha,2-\alpha),~0<\alpha\le 2&(3.26)\cr}
$$where the Riemann-Liouville operator of order $\mu$ defined by
(A3) has the initial conditions
$$^{RL}_0D_t^{\mu-1}N(x,0)=N_0(x);~
^{RL}_0D_t^{\mu-2}N(x,0)=0,0<\mu\le
2,\lim_{|x|\to\pm\infty}N(x,t)=0.\eqno(3.27)
$$and ${_xD}_{\theta}^{\alpha}$ is the Riesz-Feller space-fractional
derivative or order $\alpha$ and asymmetry $\theta$ defined by
(A11), $\eta$ is a diffusion constant and $\phi(x,t)$ is a nonlinear
function. Then for the solution of (3.25) subject to the above
constraints, there holds the formula
$$\eqalignno{N(x,t)&=\int_0^xG_5(x-\tau,t)N_0(\tau){\rm d}\tau\cr
&+\int_0^t(t-\xi)^{\beta-1}\int_0^xG_6(x-\tau,t-\xi)\phi(\tau,\xi){\rm
d}\tau~{\rm d}\xi&(3.28)\cr \noalign{\hbox{where}}
G_5(x,t)&={{t^{\mu-1}}\over{\alpha|x|}}H_{3,3}^{2,1}\left[{{|x|}\over{\eta^{{1}\over{\alpha}}t^{{\mu}\over{\alpha}}}}
\bigg\vert_{(1,{{1}\over{\alpha}}),(1,1),(1,\rho)}^{(1,{{1}\over{\alpha}}),(\mu,{{\mu}\over{\alpha}}),(1,\rho)}\right],~\alpha>0&(3.29)\cr
\noalign{\hbox{and}}
G_6(x,t)&={{1}\over{\alpha|x|}}H_{3,3}^{2,1}\left[{{|x|}\over{\eta^{{1}\over{\alpha}}t^{{\mu}\over{\alpha}}}}
\bigg\vert_{(1,{{1}\over{\alpha}}),(1,1),(1,\rho)}^{(1,{{1}\over{\alpha}}),(\mu,{{\mu}\over{\alpha}}),(1,\rho)}\right]&(3.30)\cr}
$$for $\alpha>0,\Re(\mu)>0,\rho+{{\alpha-\theta}\over{2\alpha}}$.}

\vskip.3cm\noindent{\bf 4.\hskip.3cm Fractional order moments}
\vskip.3cm In this section we will calculate the fractional order
moments defined by
$$\eqalignno{<|x(t)|^{\delta}>&=\int_{-\infty}^{\infty}|x|^{\delta}N(x,t){\rm
d}x&(4.1)\cr \noalign{\hbox{Using the result (3.9) and the following
definition of the Mellin transform}}
M\{f(t);s\}&=\int_0^{\infty}t^{s-1}f(t){\rm d}t&(4.2)\cr
\noalign{\hbox{we find}}
<|x|^{\delta}>&=\int_{-\infty}^{\infty}|x|^{\delta}N(x,t){\rm
d}x&(4.3)\cr
&={{2t^{\mu+\nu(1-\mu)-1}}\over{\alpha}}\int_0^{\infty}x^{\delta-1}H_{3,3}^{2,1}\left[{{|x|}\over{\eta^{{1}\over{\alpha}}t^{{\mu}\over{\alpha}}}}
\bigg\vert_{(1,1),(1,{{1}\over{\alpha}}),(1,\rho)}^{(1,{{1}\over{\alpha}}),(\mu+\nu(1-\mu),{{\mu}\over{\alpha}}),(1,\rho)}\right]{\rm
d}x&(4.4)\cr}
$$The following formula gives the Mellin transform of the H-function
(Mathai et al., 2010):
$$\int_0^{\infty}x^{\delta-1}H_{p,q}^{m,n}\left[ax\bigg\vert_{(b_q,B_q)}^{(a_p,A_p)}\right]{\rm
d}x=a^{-\delta}{{\{\prod_{j=1}^m\Gamma(b_j+B_j\delta)\}\{\prod_{j=1}^m\Gamma(1-a_j-A_j\delta)\}}\over{\{\prod_{j=m+1}^q\Gamma(1-b_j-B_j\delta)\}\{\prod_{j=n+1}^p\Gamma(a_j+A_j\delta)\}}}\eqno(4.5)
$$where $-\min_{1\le j\le
m}\Re({{b_j}\over{B_j}})<\Re(\delta)<\max_{1\le j\le
n}\Re({{1-a_j}\over{A_j}})$, $|\arg a|<{1\over2}\pi\theta$,
$\theta=\sum_{j=1}^mB_j-\sum_{j=m+1}^qB_j+\sum_{j=1}^nA_j-\sum_{j=n+1}^pA_j>0$.
Applying the above formula to evaluate the integral in (4.5), we see
that
$$<|x(t)|^{\delta}>={{2}\over{\alpha}}\eta^{{\delta}\over{\alpha}}t^{\mu+\nu(1-\mu)+\mu({{\delta}\over{\alpha}})-1}
{{\Gamma(-{{\delta}\over{\alpha}})\Gamma(1+\delta)\Gamma(1+{{\delta}\over{\alpha}})}\over{\Gamma(-\rho\delta)\Gamma(\mu+\nu(1-\mu)
+{{\mu\delta}\over{\alpha}})\Gamma(1+\rho\delta)}}\eqno(4.6)
$$for $-\min\{\alpha,1\}<\Re(\delta)<0$.

\vskip.3cm\noindent{\bf 5.\hskip.3cm Computational representations
of the solution (3.9)}

\vskip.3cm In this section we will derive the computational
representation of the fundamental solution (3.9), which can be
expressed in terms of the Mellin-Barnes type integral as
$$N(x,t)={{1}\over{\pi x}}{{t^{\mu+\nu(1-\mu)-1}}\over{2\pi
i}}\int_L{{\Gamma(s)\Gamma(1-s)\Gamma(1-s\alpha)}\over{\Gamma(\mu+\nu(1-\mu)-s\mu)}}
\sin[{{s\pi}\over{2\alpha}}(\theta-\alpha)][{{x^{\alpha}}\over{\eta
t^{\mu}}}]^s{\rm d}s.\eqno(5.1)
$$Let us assume that the poles of the gamma functions in the
integrand of (5.1) are all simple. Now, evaluating the sum of
residues in ascending powers of $x$ by calculating the residues at
the poles of $\Gamma(1-s)$ at the points $s=1+n, n\in N_0$ and
$\Gamma(1-s\alpha)$ at the points $s=(1+n)/\alpha,~n\in N_0$ we
obtain the following representation of the fundamental solution
(3.12) in terms of two convergent series in ascending powers of $x$
$$\eqalignno{N(x,t)&={{x^{\alpha-1}t^{\mu+\nu(1-\mu)-1}}\over{\pi(\eta
t^{\mu})}}\sum_{n=0}^{\infty}{{\Gamma[1-\alpha(1+n)]}\over{\Gamma[\mu+\nu(1-\mu)-(1+n)\mu]}}
\sin[{{(1+n)}\over{2\alpha}}(\theta-\alpha)][{{-x^{\alpha}}\over{\eta
t^{\mu}}}]^n\cr
&+{{t^{\mu+\nu(1-\mu)-1}}\over{\pi\alpha\eta^{{1}\over{\alpha}}t^{{\mu}\over{\alpha}}}}
\sum_{n=0}^{\infty}{{\Gamma[(1+n)/\alpha]\Gamma[1-(1+n)/\alpha]}\over{n!\Gamma[\mu+\nu(1-\mu)-(1+n)(\mu/\alpha)]}}
\sin[{{(1+n)\pi}\over{2\alpha^2}}(\theta-\alpha)][{{-x}\over{\eta^{{1}\over{\alpha}}t^{{\mu}\over{\alpha}}}}]^n&(5.2)\cr}
$$where $\left\vert{{x^{\alpha}}\over{\eta t^{\mu}}}\right\vert <1$.
\vskip.2cm For $\nu=1$ the time fractional derivative
$D_{0_{+}}^{\mu,\nu}$ becomes Caputo derivative $^C_0D_t^{\mu}$ and
we obtain the following simplified form of the result given by
Mainardi et al. (2001):
$$\eqalignno{N(x,t)&={{x^{\alpha-1}}\over{\pi(\eta
t^{\mu})}}\sum_{n=0}^{\infty}{{\Gamma[1-\alpha(1+n)]}\over{\Gamma[1-\mu(1+n)]}}
\sin[{{(1+n)\pi}\over{2\alpha}}(\theta-\alpha)][{{-x^{\alpha}}\over{\eta
t^{\mu}}}]^n\cr
&+{{1}\over{\pi\alpha(\eta^{{1}\over{\alpha}}t^{{\mu}\over{\alpha}})}}
\sum_{n=0}^{\infty}{{\Gamma[(1+n)/\alpha]}\over{n!\Gamma[1-(1+n)(\mu/\alpha)]}}
\sin[{{(1+n)\pi}\over{2\alpha^2}}(\theta-\alpha)][{{-x}\over{\eta^{{1}\over{\alpha}}t^{{\mu}\over{\alpha}}}}]^n&(5.3)\cr}
$$
\vskip.2cm \noindent Further, if we calculate the residues at the
poles of $\Gamma(s)$ at the points $s=-n,~n\in N_0$ it gives

$$N(x,t)={{t^{\mu+\nu(1-\mu)-1}}\over{\pi
x}}\sum_{n=0}^{\infty}{{\Gamma(1+\alpha
n)}\over{\Gamma(\mu+\nu(1-\mu)+n\mu)}}\sin[{{-n\pi}\over{2\alpha}}(\theta-\alpha)][{{-x^{\alpha}}\over{\eta~t^{\mu}}}]^{-n},~1<x<\infty.\eqno(5.4)
$$For $\nu=1$ the above result simplifies to the following series
representation:
$$N(x,t)={{1}\over{x\pi}}\sum_{n=0}^{\infty}{{\Gamma(1+\alpha
n)}\over{\Gamma(1+\mu
n)}}\sin[{{-n\pi}\over{2\alpha}}(\theta-\alpha)][-{{x^{\alpha}}\over{\eta~t^{\mu}}}]^{-n},~1<x<\infty.\eqno(5.5)
$$Finally, from (5.5) it follows that $N(x,t)\sim {{1}\over{|x|}}$
for large $|x|$.

\vskip.3cm\noindent{\bf 6.\hskip.3cm Finite number of Riesz-Feller
derivatives} \vskip.3cm Following a similar procedure it is not
difficult to establish the following

\vskip.3cm\noindent{\bf Theorem 6.1.}\hskip.3cm{\it Consider an
unified fractional reaction-diffusion model
$${_0D}_t^{\mu,\nu}N(x,t)=\sum_{j=1}^m\eta_j~{_xD}_{\theta_j}^{\alpha_j}N(x,t)+\phi(x,t)\eqno(6.1)
$$where $\eta_j,t>0,x\in R; \alpha_j,\theta_j, j=1,...,m, \mu,\nu$
are real parameters with the constraints $0<\alpha_j\le
2,~|\theta_j|\le \min\{\alpha_j,2-\alpha_j\}$, ${_0D}_t^{\mu,\nu}$
is the generalized Riemann-Liouville fractional derivative operator
defined by (A5) with the initial conditions
$$[I_{0_{+}}^{(1-\nu)(1-\mu),0}]N(x,0_{+})=N_0(x);~0<\mu<1,~0\le
\nu\le 1\eqno(6.2)
$$involving the Riemann-Liouville fractional integral of order
$(1-\nu)(1-\mu)$ evaluated for $t\to 0_{+}$ and
$\lim_{|x|\to\infty}N(x,t)=0$. Here
${_xD}_{\theta_j}^{\alpha_j},j=1,...,m$ are the Riesz-Feller space
fractional derivatives of orders $\alpha_j,j=1,...,m$ and asymmetry
$\theta_j, j=1,...,m$ respectively, defined by (A10), $\phi(x,t)$ is
a nonlinear function. Then for the solution of (6.1), subject to the
constraints above, there holds the formula
$$\eqalignno{N(x,t)&={{t^{\mu+\nu(1-\mu)-1}}\over{2\pi}}\int_{-\infty}^{\infty}N_0^{*}(k)E_{\mu,\mu+\nu(1-\mu)}(-t^{\mu}\sum_{j=1}^m\eta~\psi_{\alpha_j}^{\theta_j}(k))\exp(-ikx){\rm
d}k\cr
&+{{1}\over{2\pi}}\int_0^t\xi^{\mu-1}\int_{-\infty}^{\infty}\phi^{*}(k,t-\xi)E_{\mu,\mu}(-t^{\mu}\sum_{j=1}^m\eta_j~\psi_{\alpha_j}^{\theta_j}(k))\exp(-ikx){\rm
d}k~{\rm d}\xi.&(6.3)\cr}
$$}Some special cases of Theorem 6.1 are deduced below:
\vskip.2cm\noindent(i):~If we set $\theta_1=...=\theta_m=0$ then by
virtue of the identity (A13) we arrive at the following corollary
associated with Riesz space fractional derivative:

\vskip.3cm\noindent{\bf Corollary 6.1.}\hskip.3cm{\it Consider the
extended reaction-diffusion model
$${_0D}_t^{\mu,\nu}N(x,t)=\sum_{j=1}^m\eta_j~{_xD}_0^{\alpha_j}N(x,t)+\phi(x,t)\eqno(6.4)
$$where $\eta_j>0,j=1,...,m,t>0,x\in R, \alpha_j,j=1,...,m, \mu,\nu$ are real
parameters with the constraints $0<\alpha_j\le 2$,
$|\theta_j|\le\min(\alpha_j,2-\alpha_j)$. ${_0D}_t^{\mu,\nu}$ is the
generalized Riemann-Liouville fractional derivative operator defined
by (A7) with the initial conditions
$$[I_{0_{+}}^{(1-\nu)(1-\mu),0}]N(x,0_{+})=N_0(x); ~0<\mu<1,~0\le
\nu\le 1\eqno(6.5)
$$involving the Riemann-Liouville fractional integral of order
$(1-\nu)(1-\mu)$ evaluated for $t\to 0_{+}$ and
$\lim_{|x|\to\infty}N(x,t)=0$. Then for the solution of (6.4) there
holds the formula (6.3) with $\psi_{\alpha_j}^{\theta_j}(k)$
replaced by $|k|^{\alpha_j},j=1,...,m$.}

\vskip.3cm\noindent(ii):~If we further take $\nu=1$ in the above
Corollary 6.1 then the operator $[I_{0_{+}}^{(1-\nu)(1-\mu),0}]$
reduces to Caputo operator $^C_0D_t^{\mu}$ defined in (A5) and we
arrive at the following result:

\vskip.3cm\noindent{\bf Corollary 6.2.}\hskip.3cm{\it Consider the
extended reaction-diffusion model
$$^C_0D_t^{\mu}N(x,t)=\sum_{j=1}^m\eta_j~{_xD}_0^{\alpha_j}N(x,t)+\phi(x,t)\eqno(6.6)
$$where all the quantities are as defined above, with the initial
condition $N(x,0_{+})=N_0(x);~0<\alpha_j\le 2,~0<\mu<1$,
$\lim_{|x|\to\infty}N(x,t)=0$. ${_xD}_0^{\alpha_j},j=1,...,m$ are
the Riesz space fractional derivatives of order $\alpha_j,j=1,...,m$
defined by (A14), $\phi(x,t)$ is a nonlinear function. Then for the
solution of (6.6) there holds the formula
$$\eqalignno{
N(x,t)&={{1}\over{2\pi}}\int_{-\infty}^{\infty}N_0^{*}(k)E_{\mu,1}(-t^{\mu}\sum_{j=1}^m\eta_j~|k|^{\alpha_j})\exp(-ikx){\rm
d}k\cr
&+{{1}\over{2\pi}}\int_0^t\xi^{\mu-1}\int_{-\infty}^{\infty}\phi^{*}(k,t-\xi)E_{\mu,\mu}(-t^{\mu}\sum_{j=1}^m\eta_j~|k|^{\alpha_j})\exp(-ikx){\rm
d}k~{\rm d}\xi&(6.7)\cr}
$$}
\vskip.3cm For $m=1$ the result (6.7) reduces to one given by
Tomovski et al. (2011).

\vskip.2cm\noindent(iii):~If we set $\nu=0$ then the Hilfer (2000)
fractional operator defined by (A7) reduces to Riemann-Liouville
operator defined by (A3) and we arrive at the following

\vskip.3cm\noindent{\bf Corollary 6.3.}\hskip.3cm{\it Consider an
extended fractional reaction-diffusion model
$$^{RL}_0D_t^{\mu}N(x,t)=\sum_{j=1}^m\eta_j~{_xD}_{\theta_j}^{\alpha_j}N(x,t)+\phi(x,t)\eqno(6.8)
$$where the parameters and restrictions as defined before and with
the initial conditions
$$^{RL}_0D_t^{\mu-1}N(x,0)=N_0(x);~~
^{RL}_0D_t^{\mu-2}N(x,0)=0,~0<\mu\le
2,\lim_{|x|\to\infty}N(x,t)=0.\eqno(6.9)
$$Then for the solution of (6.8) there holds the formula:
$$\eqalignno{N(x,t)&={{t^{\mu-1}}\over{2\pi}}\int_{-\infty}^{\infty}N_0^{*}(k)E_{\mu,\mu}(-t^{\mu}\sum_{j=1}^m\eta_j~\psi_{\alpha_j}^{\theta_j}(k))\exp(-ikx){\rm
d}k\cr
&+{{1}\over{2\pi}}\int_0^t\xi^{\mu-1}\int_{-\infty}^{\infty}\phi^{*}(k,t-\xi)E_{\mu.\mu}(-t^{\mu}\sum_{j=1}^m\eta_j~\psi_{\alpha_j}^{\theta_j})\exp(-ikx){\rm
d}k~{\rm d}\xi.&(6.10)\cr}
$$}

\vskip.2cm\noindent(iv):~Finally, if we set $\theta_j=0,j=1,...,m$
in Corollary 6.3 then the Riesz-Feller derivatives reduce to Riesz
space fractional derivative and we arrive at the following

\vskip.3cm\noindent{\bf Corollary 6.4.}\hskip.3cm{\it Consider an
extended fractional reaction-diffusion model
$$^{RL}_0D_t^{\mu}N(x,t)=\sum_{j=1}^m\eta_j~{_xD}_0^{\alpha_j}N(x,t)+\phi(x,t)\eqno(6.11)
$$with the parameters and conditions on them as defined before and
with the initial conditions as in (6.7) then for the solution of
(6.11) there holds the formula
$$\eqalignno{N(x,t)&={{t^{\mu-1}}\over{2\pi}}\int_{-\infty}^{\infty}N_0^{*}(k)E_{\mu,\mu}(-t^{\mu}\sum_{j=1}^m\eta_j~|k|^{\alpha_j})\exp(-ikx){\rm
d}k\cr
&+{{1}\over{2\pi}}\int_0^t\xi^{\mu-1}\int_{-\infty}^{\infty}\phi^{*}(k,t-\xi)E_{\mu.\mu}(-t^{\mu}\sum_{j=1}^m\eta_j~|k|^{\alpha_j})\exp(-ikx){\rm
d}k~{\rm d}\xi.&(6.13)\cr}
$$}

\vskip.3cm\noindent{\bf 7.\hskip.3cm Conclusions}

\vskip.3cm In this paper, the authors have presented an extension of
the fundamental solution of space-time fractional diffusion given by
Mainardi-Luchko-Pagnini (2001) by using the fractional order
derivative operator defined by Hilfer (2000). The fundamental solution of
the equation (2.1) is obtained in terms of H-function in closed and
computable forms. Computational representations and fractional
moments of the solutions are also obtained which will enhance the
utility of the derived results in practical problems. Solutions of
unified reaction-diffusion models associated with a finite number of
Riesz-Feller space fractional derivatives are also investigated.

\vskip.3cm\noindent {\bf Acknowledgment} \vskip.3cm The authors
would like to thank the Department of Science and Technology,
Government of India, for the financial support for this work under
project No. SR/S4/MS:287/05.

\vskip.3cm\noindent{\bf References}

\vskip.3cm\noindent Caputo, M. (1969):~~{\it Elasticita e
Dissipazione}, Zanichelli, Bologna.

\vskip.2cm\noindent Del-Castillo-Negrete, D., Carreras, B.A. and Lynch,
V. (2003):~~Front dynamics in diffusion systems with L\'evy
flights: a fractional, diffusion approach, {\it Physical Review
Letters}, {\bf 91}, 01832.

\vskip.2cm\noindent Del-Castillo-Negrete, Carreas, B.A. and Lynch,
V. (2002):~~Front propagation and segregation in a
reaction-diffusion model with cross-diffusion, {\it Physica D}, {\bf
168-169}, 45-60.

\vskip.2cm\noindent Dzherbashyan, M.M. (1966):~~{\it Integral
Transforms and Representation of Functions in Complex Domain (in
Russian)}, Nauka, Moscow.

\vskip.2cm\noindent Dzherbashyan, M.M. (1993):~~{\it Harmonic
Analysis and Boundary Value Problems in the Complex Domain},
Birkhaeuser, Basel.

\vskip.2cm\noindent Erd\'elyi, A., Magnus, W., Oberhettinger, F. and
Tricomi, F.G. (1954):~~{\it Tables of Integral Transforms}, Vol.1,
McGraw-Hill, New York.

\vskip.2cm\noindent Erd\'elyi, A., Magnus, W., Oberhettinger, F. and
Tricomi, F.G. (1955):~~{\it Higher Transcendental Functions},
Vol.3, McGraw-Hill, New York.

\vskip.2cm\noindent Feller, W. (1952):~~On a generalization of
Marcel Riesz' potentials and the semi-groups generated by them. {\it
Meddelanden Lunds Universitets Matematiska Seminarium Comm. S\'em
Math\'em Universit\'e de Lunds (Suppl. d\'edi\'e \'a M. Riesz.
Lund)} 73-81.

\vskip.2cm\noindent Feller W. (1971):~~{\it An Introduction to
Probability and Its Applications}, Vol.2, 2nd edition, Wiley, New
York (1st edition 1966).

\vskip.2cm\noindent Gorenflo, R. and Mainardi, F.
(1999):~~Approximation of L\'evy-Feller diffusion by random walk,
{\it Journal for Analysis and its Applications}, {\bf 18(2)}, 1-16.

\vskip.2cm\noindent Gorenflo, R. and Mainardi, F.
(1997):~~Fractional calculus integral and differential equations of
fractional order, in A. Cparinteri and F. Mainardi (eidtors), {\it
Fractals and Fractional Calculus in Continuum Mechanics}, Wien and
New York, Springer, pp. 223-276.

\vskip.2cm\noindent Haken, H. (2004):~~{\it Synergetics Introduction
and Advanced Topics}, Springer, Berlin-Heidelberg.

\vskip.2cm\noindent Haubold, H.J. and Mathai, A.M. (1995):~~ A
heuristic remark on the periodic variation in the number of solar
neutrinos detected on Earth, {\it Astrophysics and Space Science},
{\bf 228}, 113-134.

\vskip.2cm\noindent Haubold, H.J. and Mathai, A.M. (2000):~~The
fractional kinetic equation and thermonuclear functions, {\it
Astrophysics and  Space Science}, {\bf 273}, 53-63.

\vskip.2cm\noindent Haubold, H.J., Mathai, A.M. and Saxena, R.K.
(2007):~~Solutions of the reaction-diffusion equations in terms of
the H-functions, {\it Bulletin Astro. Soc. India}, {\bf 35(4)},
681-689.

\vskip.2cm\noindent Haubold, H.J., Mathai, A.M. and Saxena, R.K.
(2011):~~Further solutions of reaction-diffusion equations in terms
of the H-function, {\it J. Comput. Appl. Math.}, {\bf 235},
1311-1316.

\vskip.2cm\noindent Haubold, H.J., Mathai, A.M., and Saxena, R.K. (2012):
~~Analysis of solar neutrino data from SuperKamiokande I and II: Back
to the solar neutrino problem, arXiv:astro-ph.SR/1209.1520.

\vskip.2cm\noindent Henry, B.I. and Wearne, S.L. (2000):~~Fractional
reaction-diffusion, {\it Physica A}, {\bf 276}, 448-455.

\vskip.2cm\noindent Henry, B.I. and Wearne, S.L. (2002):~~Existence
of Turing instabilities in a two-species fractional
reaction-diffusion system, {\it SIAM Journal of Applied
Mathematics}, {\bf 62}, 870-887.

\vskip.2cm\noindent Henry, B.I., Langlands, T.A.M. and Wearne, S.L.
(2005):~Turing pattern formation in fractional activator-inhibitor
systems, {\it Physical Review E}, {\bf 72}, 026101.

\vskip.2cm\noindent Hilfer, R. (2000):~~Fractional time evolution,
In: Hilfer, R. (editor), {\it Applications of Fractional Calculus in
Physics}, World Scientific Publishing, Singapore, pp. 87-130.

\vskip.2cm\noindent Hilfer, R. (2003):~~ On fractional relaxation,
{\it Fractals}, {\bf 11}, 251-257.

\vskip.2cm\noindent Hilfer, R. (2009):~~Threefold Introduction
to Fractional Derivatives, in {\it Anomalous Transport: Foundations \&
Applications}, Edited by R. Klages, G. Radons and I.M.
Sokolov, Wiley-VCH, Weinheim.

\vskip.2cm\noindent Hundsdorfer, W. and Verwer, J.G. (2003):~~{\it
Numerical Solution of Time-Dependent Advection-Diffusion-Reaction
Equations}, Springer, Berlin-Heidelberg-New York.

\vskip.2cm\noindent Jespersen, S., Metzler, R. and Fodgeby, H.C.
(1999):~~L\'evy flights in external force fields: Langevin and
fractional Fokker-Planck equations and their solutions, {\it
Physical Review E}, {\bf 59(3)}, 2736-2745.

\vskip.2cm\noindent Kilbas, A.A., Peirentozi, T., Trujillo, J.J. and
Vazquez, L. (2004):~~On the solution of fractional evolution
equation, {\it J. Phys. A, Math. Gen.,}{\bf 37(9)}, 3272-3283.

\vskip.2cm\noindent Kilbas, A.A., Srivastava, H.M. and Trujillo,J.J.
(2006):~~{\it Theory and Applications of Fractional Differential
Equations}, Elsevier, Amsterdam.

\vskip.2cm\noindent Kuramoto, Y. (2003):~~{\it Chemical
Oscillations, Waves, and Turbulence}, Dover Publications, Mineola,
New York.

\vskip.2cm\noindent Mainardi, F. , Luchko, Y. and Pagnini, G.
(2001):~~The fundamental solution of the space-time fractional
diffusion equation, {\it Fract. Calc. Appl. Anal.}, {\bf 4(2)},
153-192.

\vskip.2cm\noindent Mainardi, F., Pagnini, G. and Saxena, R.K.
(2005):~~ Fox H-functions in fractional diffusion, {\it J. Comput.
Appl. Math.}, {\bf 178}, 321-331.

\vskip.2cm\noindent Mathai, A.M. and Saxena, R.K. (1978):~~{\it The
H-function with Applications in Statistics and Other Disciplines},
Wiley Eastern New Delhi and Wiley Halsted, New York.

\vskip.2cm\noindent Mathai, A.M., Saxena, R.K. and Haubold, H.J.
(2010):~~{\it The H-function: Theory and Applications}, Springer,
New York.

\vskip.2cm\noindent Miller, K.S. and Ross, B. (1993):~~{\it An
Introduction to the Fractional Calculus and Fractional Differential
Equations}, Wiley, New York.

\vskip.2cm\noindent Mittag-Leffler, G.M. (1903):~~ Sur la nouvelle
fonction $E_{\alpha}(x)$, {\it C.R. Acad. Sci., Paris (ser.ii)},
{\bf 137}, 554-558.

\vskip.2cm\noindent Mittag-Leffler, G.M. (1905):~~Sur la
representation analytique d'une fonction branche uniforme d'une
fonction, {\it Acta Math.}, {\bf 239}, 101-181.

\vskip.2cm\noindent Murray, J.D. (2003):~~{\it Mathematical
Biology}, Springer, New York.

\vskip.2cm\noindent Nicolis, G. and Prigogine, I. (1997):~~{\it
Self-Organization in Nonequilibrium Systems: From Dissipative
Structures to Order Through Fluctuations}, Wiley, New York.

\vskip.2cm\noindent Oldham, K.B. and Spanier, J. (1974):~~{\it The
Fractional Calculus: Theory and Applications of Differentiation and
Integration of Arbitrary Order}, Academic Press, New York.

\vskip.2cm\noindent Prudnikov, A.P., Brychkov, Yu. A. and Marichev,
O.I. (1989):~~{\it Integrals and Series, Vol.3, More Special
Functions}, Gordon and Breach, New York.

\vskip.2cm\noindent Samko, S.G., Kilbas, A.A. and Marichev, O.I.
(1990):~~{\it Fractional Integrals and Derivatives: Theory and
Applications}, Gordon and Breach, New York.

\vskip.2cm\noindent Saxena, R.K. (2012):~~Solution of fractional
partial differential equation related to quantum mechanics, {\it
Algebras, Groups and Geometries}, (to appear).

\vskip.2cm\noindent Saxena, R.K., Mathai, A.M. and Haubold, H.J.
(2002):~~On fractional kinetic equations, {\it Astrophysics and
Space Science}, {\bf 282}, 281-287.

\vskip.2cm\noindent Saxena, R.K., Mathai, A.M. and Haubold, H.J.
(2004a):~~On generalized fractional kinetic equations, {\it Physica
A}, {\bf 344}, 657-664.

\vskip.2cm\noindent Saxena, R.K., Mathai, A.M. and Haubold, H.J.
(2004b):~~Unified fractional kinetic equation and a fractional
diffusion equation, {\it Astrophysics and Space Science}, {\bf
290(3\&4)}, 299-310.

\vskip.2cm\noindent Saxena, R.K., Mathai, A.M. and Haubold, H.J.
(2004c):~~Astrophysical thermonuclear functions for Boltzmann-Gibbs
statistics and Tsallis statistics, {\it Physica A}, {\bf 344},
649-656.

\vskip.2cm\noindent Saxena, R.K., Mathai, A.M. and Haubold, H.J.
(2006a):~~ Fractional reaction-diffusion equations, {\it
Astrophysics an Space Science}, {\bf 305}, 289-296.

\vskip.2cm\noindent Saxena, R.K., Mathai, A.M. and Haubold, H.J.
(2006b):~~Reaction-diffusion systems and nonlinear waves, {\it
Astrophysics and Space Science}, {\bf 305}, 297-303.

\vskip.2cm\noindent Saxena, R.K., Mathai, A.M. and Haubold, H.J.
(2006c):~~ Solution of generalized fractional reaction-diffusion
equations, {\it Astrophysics and Space Science}, {\bf 305}, 305-313.

\vskip.2cm\noindent Saxena, R.K., Mathai, A.M. and Haubold, H.J.
(2006d):~~Solutions of fractional reaction-diffusion equations in
terms of the Mittag-Leffler functions, {\it Int. J. Sci. Res.}, {\bf
15}, 1-17.

\vskip.2cm\noindent Saxena, R.K., Mathai, A.M. and Haubold, H.J.
(2008):~~ Solution of a fractional kinetic equation and a fraction
diffusion equation, {\it Int. J. Sci. Res.}, {\bf 17}, 1-8.

\vskip.2cm\noindent Saxena R.K., Saxena, R. and Kalla, S.L.
(2010):~~Solution of space-time fractional Schr\"odinger equation
occurring in quantum mechanics, {\it Fractional Calculus and Applied
Analysis}, {\bf 13(2)}, 177-190.

\vskip.2cm\noindent Schneider, W.R. (1986):~~Stable distributions:
Fox function representation and generalization, S. Albeverio, G.
Casati and D. Merilini (editors), {\it Stochastic Processes in
Classical and Quantum Systems}, Berlin-Heidelberg, Springer,
pp. 497-511 [Lecture Notes in Physics, Vol. 262].

\vskip.2cm\noindent Srivastava, H.M. and Tomovski, Z.
(2009):~~Fractional calculus with   an integral operator containing
a generalized Mittag-Leffler function in the kernel, {\it Appl.
Math. Comput.}, {\bf 211}, 198-210.

\vskip.2cm\noindent Tomovski, Z., Sandev, T., Metzler, R. and
Dubbeldam, J. (2011):~~Generalized space-time fractional diffusion
equation with composite fractional time derivatives, {\it Physica
A}, doi:10.10`6/j.physa, 2011.12.035.

\vskip.2cm\noindent Uchaikin, V.V. and Zolotarev, V.M. (1999):~~{\it
Chance and Stability: Stable Distributions \& Their Applications},
Utrecht, VSP.

\vskip.2cm\noindent Wihelmsson, H. and Lazzaro, E. (2001):~~{\it
Reaction-Diffusion Problems in the Physics of Hot Plasmas},
Institute of Physics Publishing, Bristol and Philadelphia.

\vskip.2cm\noindent Wiman, A. (1905):~~Ueber den Fundamentalsatz in
der Theorie der Functionen, {\it Acta Math.}, {\bf 29}, 191-201.

\vskip.3cm\noindent{\bf Appendix A.\hskip.3cm Mathematical
preliminaries}

\vskip.3cm A generalization of the Mittag-Leffler function
$$E_{\alpha}(z)=\sum_{k=0}^{\infty}{{z^k}\over{\Gamma(1+k\alpha)}},~\alpha\in
C, \Re(\alpha)>0\eqno(A1)
$$was introduced by Wiman (1905) in the generalized form
$$E_{\alpha,\beta}(z)=\sum_{k=0}^{\infty}{{z^k}\over{\Gamma(\beta+\alpha
k)}},~\Re(\alpha)>0, \Re(\beta)>0.\eqno(A2)
$$The main results of these functions are available in the handbook
of Erd\'elyi, et al. (1955, Section 18.1) and the monographs of
Dzherbashyan (1966, 1993). The left-sided Riemann-Liouville
fractional integral of order $\nu$ is defined by Miller and Ross
(1993, p.45), Samko et al. (1990), Kilbas et al. (2006) as
$$^{RL}_0D_t^{-\nu}N(x,t)=I_0^{\nu}N(x,t)={{1}\over{\Gamma(\nu)}}\int_0^t(t-u)^{\nu-1}N(x,u){\rm
d}u,~t>0,\Re(\nu)>0.\eqno(A3)
$$The left-sided Riemann-Liouville fractional derivative of order
$\alpha$ is defined as
$$^{RL}_0D_t^{\mu}N(x,t)=[{{{\rm d}}\over{{\rm
d}x}}]^n(I_0^{n-\mu}N(x,t)),~\Re(\mu)>0,n=[\Re(\mu)]+1\eqno(A4)
$$where $[x]$ represents the greatest integer in the real number $x$.
Caputo derivative (Caputo, 1969) is defined in the form
$$\eqalignno{^C_0D_t^{\alpha}f(x,t)&={{1}\over{\Gamma(m-\alpha)}}\int_0^t{{f^{(m)}(x,t)}\over{(t-\tau)^{\alpha+1-m}}}{\rm
d}\tau,~m-1<\Re(\alpha)<m,~m\in N&(A5)\cr
&={{\partial^mf(x,t)}\over{\partial t^m}},\hbox{ for
}\alpha=m&(A6)\cr}
$$where ${{\partial^m}\over{\partial t^m}}f(x,t)$ is the $m$-th
partial derivative of $f(x,t)$ with respect to $t$. When there is no
confusion, then the Caputo operator $^C_0D_t^{\alpha}$ will be
simply denoted by ${_0D_t^{\alpha}}$. \vskip.2cm A generalization of
the Riemann-Liouville fractional derivative operator (A4) as well as
Caputo fractional derivative operator (A5) is given by Hilfer (2000)
by introducing a left-sided fractional derivative operator of two
parameters of order $0< \mu<1$ and type $0\le\nu\le 1$ in the form
$${_0D}_{a_{+}}^{\mu,\nu}N(x,t)=\left[I_{a_{+}}^{\nu(1-\mu)}{{\partial}\over{\partial
x}}\left(I_{a_{+}}^{(1-\nu)(1-\mu)}N(x,t)\right)\right].\eqno(A7)
$$For $\nu=0$, (A7) reduces to the classical Riemann-Liouville fractional derivative operator. On the other hand
for $\nu=1$ it yields the Caputo fractional derivative operator
defined by (A5). The Laplace transform formula for this operator is
given by Hilfer (2000).
$$L[{_0D}_{0_{+}}^{\mu,\nu}N(x,t);s]=s^{\mu}\tilde{N}(x,s)-s^{\nu(\mu-1)}I_{0_{+}}^{(1-\nu)(1-\mu)}N(x,0_{+}),~0<\mu<1\eqno(A8)
$$where the initial value term
$I_{0_{+}}^{(1-\nu)(1-\mu)}N(x,0_{+})$ involves the
Riemann-Liouville fractional integral operator of order
$(1-\nu)(1-\mu)$ evaluated in the limit as $t\to 0_{+}$, it being
understood that the integral
$$\tilde{N}(x,s)=L\{N(x,t);s\}=\int_0^{\infty}{\rm
e}^{-st}N(x,t){\rm d}t,\eqno(A9)
$$ where $\Re(s)>0$, exists.

\vskip.3cm\noindent{\bf Note A1.}\hskip.3cm The derivative defined
by (A7) also occurs in recent papers by Hilfer (2003, 2009),
Srivastava et al.(2009), Tomovski et al. (2011) and Saxena et al.
(2010). \vskip.2cm Following Feller (1952, 1971), it is conventional
to define the Riesz-Feller space fractional derivative of order
$\alpha$ and skewness $\theta$ in terms of its Fourier transform as
$$\eqalignno{F\{{_xD}_{\theta}^{\alpha}f(x);k\}&=-\psi_{\alpha}^{\theta}(k)f^{*}(k),&(A10)\cr
\noalign{\hbox{where}}
\psi_{\alpha}^{\theta}(k)&=|k|^{\alpha}\exp[i(sign~k){{\theta\pi}\over{2}}],~0<\alpha\le
2. ~|\theta|\le \min\{\alpha,2-\alpha\}.&(A11)\cr}
$$When $\theta=0$ we have a symmetric operator with respect to $x$,
that can be interpreted as
$${_xD}_{\theta}^{\alpha}=-[-{{{\rm d}^2}\over{{\rm
d}x^2}}]^{{\alpha}\over2}.\eqno(A12)
$$This can be formally deduced by writing
$-(k)^{\alpha}=-(k^2)^{{\alpha}\over2}$. For $\theta=0$ we also have
$$F\{{_xD}_0^{\alpha}f(x);k\}=-|k|^{\alpha}f^{*}(k).\eqno(A13)
$$For $0<\alpha\le 2$ and $|\theta|\le\min\{\alpha,2-\alpha\}$ the
Riesz-Feller derivative can be shown to possess the following
integral representation in $x$ domain:
$$\eqalignno{{_xD}_{\theta}^{\alpha}f(x)&={{\Gamma(1+\alpha)}\over{\pi}}\bigg\{\sin[(\alpha+\theta){{\pi}\over2}]\int_0^{\infty}{{f(x+\xi)-f(x)}\over{\xi^{1+\alpha}}}{\rm
d}\xi\cr
&+\sin[(\alpha-\theta){{\pi}\over2}]\int_0^{\infty}{{f(x-\xi)-f(x)}\over{\xi^{1+\alpha}}}{\rm
d}\xi\bigg\}.\cr}
$$For $\theta=0$, the Riesz-Feller fractional derivative becomes the
Riesz fractional derivative of order $\alpha$ for $1<\alpha\le 2$
defined by analytic continuation in the whole range $0<\alpha\le 2$,
$\alpha\ne 1$, see Gorenflo and Mainardi (1999), as
$$\eqalignno{{_xD}_{0}^{\alpha}&=-\lambda[I_{+}^{-\alpha}-I_{-}^{-\alpha}]&(A14)\cr
\noalign{\hbox{where}} \lambda&={{1}\over{2\cos(\alpha\pi/2)}};
~~I_{\pm}^{-\alpha}={{{\rm d}^2}\over{{\rm
d}x^2}}I_{\pm}^{2-\alpha}.&(A15)\cr}
$$The Weyl fractional integral operators are defined in the
monograph by Samko et al. (1990) as
$$\eqalignno{(I_{+}^{\beta}
N)(x)&={{1}\over{\Gamma(\beta)}}\int_{-\infty}^x(x-\xi)^{\beta-1}N(\xi){\rm
d}\xi,~\beta>0;&(A16)\cr
(I_{-}^{\beta}N)(x)&={{1}\over{\Gamma(\beta)}}\int_x^{\infty}(\xi-x)^{\beta-1}N(\xi){\rm
d}\xi,~\beta>0.&(A17)\cr}
$$
\vskip.3cm\noindent{\bf Note A2}.\hskip.3cm We note that
${_xD}_0^{\alpha}$ is a pseudo differential operator. In particular,
we have
$${_xD}_0^2={{{\rm d}^2}\over{{\rm d}x^2}},\hbox{  but
}{_xD}_0^1\ne {{{\rm d}}\over{{\rm d}x}}.\eqno(A18)
$$

\bye